\documentclass[11pt,reqno,tbtags]{amsart}
\usepackage{amssymb}
\usepackage{natbib}

\usepackage{graphicx}
\usepackage{psfrag}

\bibpunct[, ]{[}{]}{;}{n}{,}{,}
\numberwithin{equation}{section}
\allowdisplaybreaks

\newtheorem{thm}{Theorem}[section]
\newtheorem{lem}[thm]{Lemma}

\newtheorem{cor}[thm]{Corollary}
\newtheorem{conj}[thm]{Conjecture}

\theoremstyle{definition}
\newtheorem{rem}[thm]{Remark}

\theoremstyle{remark}

\newenvironment{romenumerate}{\begin{enumerate}% gives (i), (ii) etc.
 }{\end{enumerate}}

\newcounter{oldenumi}
{\setcounter{oldenumi}{\value{enumi}}
\begin{romenumerate} \setcounter{enumi}{\value{oldenumi}}}
{\end{romenumerate}}

\newcounter{thmenumerate}

\newcounter{xenumerate}   %no left indentation; thus wider lines

\newcommand{\refT}[1]{Theorem~\ref{#1}}

\newcommand{\refL}[1]{Lemma~\ref{#1}}

\newcommand{\refS}[1]{Section~\ref{#1}}

\newcommand{\refF}[1]{Figure~\ref{#1}}

\newcommand\marginal[1]{\marginpar{\raggedright\parindent=0pt\tiny #1}}

\begingroup
  \divide\count255 by 60
  \multiply\count255 by -60
  \advance\count255 by \time
  \ifnum \count255 < 10 \xdef\klockan{\the\count1.0\the\count255}
  \else\xdef\klockan{\the\count1.\the\count255}\fi
\endgroup

\newcommand\set[1]{\ensuremath{\{#1\}}}
\newcommand\bigset[1]{\ensuremath{\bigl\{#1\bigr\}}}
\newcommand\Bigset[1]{\ensuremath{\Bigl\{#1\Bigr\}}}

\newcommand\bigpar[1]{\bigl(#1\bigr)}
\newcommand\Bigpar[1]{\Bigl(#1\Bigr)}

\newcommand\lrpar[1]{\left(#1\right)}

\def\rompar(#1){\textup(#1\textup)}    % usage: \rompar(...)

\newcommand\parfrac[2]{\Bigpar{\frac{#1}{#2}}}

\newcommand\expe[1]{e^{#1}}
\def\xexp(#1){e^{#1}}

\newcommand\iid{i.i.d.\spacefactor=1000}    
\newcommand\ie{i.e.\spacefactor=1000}
\newcommand\eg{e.g.\spacefactor=1000}

\newcommand\cf{cf.\spacefactor=1000}
\newcommand{\as}{a.s.\spacefactor=1000}

\newcommand\eqd{\overset{\mathrm{d}}{=}}

\newcommand\bbN{\mathbb N}

\newcounter{CC} 
\newcounter{cc}
\newcommand\E{\operatorname{\mathbb E{}}}
\renewcommand\P{\operatorname{\mathbb P{}}}

\newcommand\Exp{\operatorname{Exp}}

\newcommand\ga{\alpha}

\newcommand\gG{\Gamma}

\newcommand\gl{\lambda}
\newcommand\gL{\Lambda}

\newcommand\cA{\mathcal A}

\newcommand\fS{\mathfrak S}

\newcommand\tgl{{\widetilde \gl}}
\newcommand\tXi{{\widetilde \Xi}}

\def\[#1]{[\![#1]\!]}

\newcommand\qq{^{1/2}}
\newcommand\qqw{^{-1/2}}
\newcommand\qw{^{-1}}

\newcommand\qqc{^{3/2}}
\newcommand\qqcw{^{-3/2}}

\renewcommand{\=}{:=}

\newcommand\intoi{\int_0^1}

\newcommand\intotau{\int_0^T}
\newcommand\intoo{\int_0^\infty}

\newcommand\oi{[0,1]}
\newcommand\ooo{[0,\infty)}

\newcommand\dd{\,\textup{d}}

\newcommand\tpo{_{t\ge0}}
\newcommand\nnu{_{\nu}}
\newcommand\nnux{_{\nu^*}}
\newcommand\inu{I\nnu}
\newcommand\exc{\mathbf e}
\newcommand\enu{\exc\nnu}
\newcommand\enux{\exc\nnux}
\newcommand\nuioo{_{\nu=1}^\infty}
\newcommand\iioo{_{i=1}^\infty}
\newcommand\tx{\zeta}
\newcommand\taux{T}
\newcommand\fa{f_{\cA}}

\newcommand\REM[1]{{\raggedright\texttt{[#1]}\par\marginal{XXX}}}

\newcommand\urladdrx[1]{{\urladdr{\def~{{\tiny$\sim$}}#1}}}

\begin{document}
\title%[]
{The integral of the supremum process of Brownian motion}

\date{June 28, 2007} % (typeset \today{} \klockan)} %October 19, 2004}

\author{Svante Janson}
\address{Department of Mathematics, Uppsala University, PO Box 480,
SE-751~06 Uppsala, Sweden}
\email{svante.janson@math.uu.se}
\urladdrx{http://www.math.uu.se/~svante/}
\author{Niclas Petersson}
\address{Department of Mathematics, Uppsala University, PO Box 480,
SE-751~06 Uppsala, Sweden}
\email{niclas.petersson@math.uu.se}
\urladdrx{http://www.math.uu.se/~niclasp/}

%{60C05 (68P10,68W40)} %%{Primary: <subject>; Secondary: <subject>}

\begin{abstract}
In this paper we study the integral of the supremum
process of standard Brownian motion. We present an explicit formula
for the moments of the integral (or area) $\cA(T)$, covered by the
process in the time interval $[0,T]$. The Laplace transform of
$\cA(T)$ follows as a consequence. The main proof involves a
double Laplace transform of $\cA(T)$ and is based on excursion
theory and local time for Brownian motion. 
\end{abstract}

\maketitle

\section{Introduction}\label{S:intro}

Let $B(t)$, $t\ge0$, be a standard Brownian motion.
Consider the following associated processes: the supremum process
$S(t)\=\max_{0\le s\le t} B(t)$, and the local time $L(t)$, which can
be regarded as a measure of the time $B(t)$ spends at 0 in the interval
$[0,t]$, see Revuz and Yor \cite[Chapter VI]{RY} for details. It is well-known that these two processes, although pathwise quite
different, have the same
distribution \cite[Chapter VI.2]{RY},
\begin{equation*}
\bigset{S(t)}_{t\ge0} \eqd \bigset{L(t)}_{t\ge0}.
\end{equation*}

The purpose of this paper is to study the distribution of the 
area under $S(t)$ or, equivalently, $L(t)$ over a given time interval
$[0,\taux]$. That is, the integral
\begin{equation}\label{a1}
  \cA(\taux)\=\intotau S(t)\dd t \eqd
\intotau L(t)\dd t.
\end{equation}
For ease of notation, let $\cA\=\cA(1)$.

The area (\ref{a1}) appeared as a random parameter when analysing
displacements for linear probing hashing. The Laplace transform of
$\cA$, which is presented in Corollary \ref{C1}, provided the means to
prove one of the main theorems in Petersson \cite{NP03}.  

Note that the usual Brownian scaling
\begin{equation*}
\bigset{B(\taux t)}\tpo\eqd\bigset{\taux\qq B(t)}\tpo,
\end{equation*}
for any $T>0$, implies the corresponding scaling for the supremum process,
\begin{equation*}
\bigset{S(\taux t)}\tpo\eqd\bigset{\taux\qq S(t)}\tpo.
\end{equation*}
Thus, for $T>0$,
\begin{equation}\label{a2}
\cA(\taux) =  \taux \intoi S(\taux t)\dd t
\eqd \taux\qqc \cA,
\end{equation}
and it is enough to study $\cA$.

\section{Results}\label{S:res}

Let $\psi(s)\=\E e^ {-s\cA}$ denote the Laplace transform of
$\cA$. 
An essential part of this paper is devoted to proving the
following formula for the Laplace transform of a variation of
$\psi$, 
or in other words, a \emph{double} Laplace transform of
$\cA$. Such formulas have already been derived for the integral of
$|B(t)|$ and other similar integrals of processes related to Brownian
motion, see Perman and Wellner \cite{PW} and the survey by Janson
\cite{SJ201}.  

\begin{thm}\label{T1}
Let $\psi$ be the Laplace transform of $\cA$. For all $\ga,\gl>0$,
\begin{equation*}
\intoo \psi\bigpar{\ga s\qqc} e^{-\gl s} \dd s 
= \intoo\Bigpar{1+\frac{3 \ga s}{2 \sqrt{2 \gl}}}^{-2/3} e^{-\gl s} \dd s. 
\end{equation*}
\end{thm}

\begin{rem}
One of the parameters $\ga$ and $\gl$ in \refT{T1} can be eliminated
(by setting it equal to 1, for instance)
without loss of generality. In fact, for any $\beta>0$, 
the formula is preserved by the substitutions
$\gl\mapsto \beta \gl$, $\ga\mapsto \beta \qqc\ga$ and $s\mapsto \beta \qw s$.
\end{rem}

The proof is given in \refS{Spf}. It is based on excursion theory for
Brownian motion and is inspired by similar arguments for other
Brownian areas, see Perman and Wellner \cite{PW}.

\begin{thm} \label{T2}
The n:th moment of $\cA$ is
\begin{equation*}
\E \cA^n = \frac{n! \, \gG(n+2/3)}{\gG(2/3) \, \gG(3n/2+1)}
\Bigpar{\frac{3 \sqrt{2}}{4}}^n, \qquad n \in \bbN. 
\end{equation*}
\end{thm}

\begin{proof} Set $\gl = 1$ in \refT{T1} 
and denote the left and right hand side by
\begin{align*}
I(\ga) & := \intoo \psi\bigpar{\ga s^{3/2}} e^{-s} \dd s
\intertext{and}
J(\ga) & := \intoo \Bigpar{1+\frac{3 \ga s}{2 \sqrt{2}}}^{-2/3} e^{-s} \dd s.
\end{align*}

The integrand of $I(\ga)$ and all its derivatives with respect to
$\ga$ are dominated by functions of the form
$s^Ke^{-s}$, uniformly in $\ga>0$. Differentiation of 
$I(\ga)$ is therefore allowed indefinitely due to dominated
convergence. The same argument applies to $J(\ga)$.  

Also, the dominated convergence theorem shows that integration 
(with respect to $s$) can be interchanged with taking the limit $\ga \to 0+$. 
Thus
\begin{align*}
\lim_{\ga \rightarrow 0+} \frac{\dd^n I(\ga)}{\dd\ga^n} & = \lim_{\ga \rightarrow 0+} \intoo \frac{\dd^n}{\dd\ga^n} \psi\bigpar{\ga s^{3/2}} e^{-s} \dd s\\
& = \intoo  \lim_{\ga \rightarrow 0+} (-s^{3/2})^n \E \Bigpar{ \cA^n \exp \bigset{-\ga s^{3/2} \cA }} e^{-s} \dd s\\
& = (-1)^n \E \bigpar{\cA^n} \intoo s^{3n/2}e^{-s} \dd s \\
& = (-1)^n \, \gG(3n/2+1) \E \cA^n
\intertext{and}
\lim_{\ga \rightarrow 0+} \frac{\dd^n J(\ga)}{\dd\ga^n} & = \lim_{\ga \rightarrow 0+} \intoo \frac{\dd^n}{\dd\ga^n} \Bigpar{1+\frac{3 \ga s}{2 \sqrt{2}}}^{-2/3} e^{-s} \dd s \\
& = \intoo  \lim_{\ga \rightarrow 0+} \frac{\gG(n+2/3)}{\gG(2/3)} \Bigpar{\frac{-3 s}{2 \sqrt{2}}}^n \Bigpar{1+\frac{3 \ga s}{2 \sqrt{2}}}^{-n-2/3} e^{-s} \dd s \\
& = \frac{\gG(n+2/3)}{\gG(2/3)} \Bigpar{\frac{-3}{2 \sqrt{2}}}^n \intoo s^{n} e^{-s} \dd s \\ 
& = \frac{\gG(n+2/3)}{\gG(2/3)} \Bigpar{\frac{-3 \sqrt{2}}{4}}^n  n!. 
\end{align*}
The fact that $I(\ga)=J(\ga)$ completes the proof.
\end{proof}

The first four moments of $\cA$ are listed in Table \ref{Tabmom}. Further, Stirling's formula provides the asymptotic relation \begin{equation}\label{moms}
\E \cA^n \sim \frac{2 \sqrt{3 \pi} }{3 \, \gG(2/3)} n^{1/6}
\Bigpar{\frac{n}{3 e}}^{n/2}, \qquad  n \to \infty.
\end{equation}

\begin{cor}\label{C1} 
The Laplace transform of $\cA$ is
\begin{equation}
\psi(s) 
= \frac{1}{\gG(2/3)} \sum_{n=0}^{\infty}
\frac{\gG(n+2/3)}{\gG(3n/2+1)} \Bigpar{\frac{-3 \sqrt{2} \, s}{4}}^n. 
\end{equation}
\end{cor}
\begin{proof} The corollary follows from the identity
\begin{equation*}
\psi(s) = \sum_{n=0}^{\infty} \frac{(-s)^n}{n!} \E \cA^n.
\end{equation*}
Note that the sum converges absolutely for every complex $s$.
\end{proof}

\begin{figure}[ht]
\psfrag{y}{$\psi(s)$}
\psfrag{x}{\hspace{1.4mm} $s$}
\psfrag{1,0}{{\scriptsize$1.0$}}
\psfrag{0,8}{{\scriptsize$0.8$}}
\psfrag{0,6}{{\scriptsize$0.6$}}
\psfrag{0,4}{{\scriptsize$0.4$}}
\psfrag{0,2}{{\scriptsize$0.2$}}
\psfrag{0}{{\scriptsize$0$}}
\psfrag{2}{{\scriptsize$2$}}
\psfrag{4}{{\scriptsize$4$}}
\psfrag{6}{{\scriptsize$6$}}
\psfrag{8}{{\scriptsize$8$}}
\begin{center}
\includegraphics[width=10cm,height=5cm]{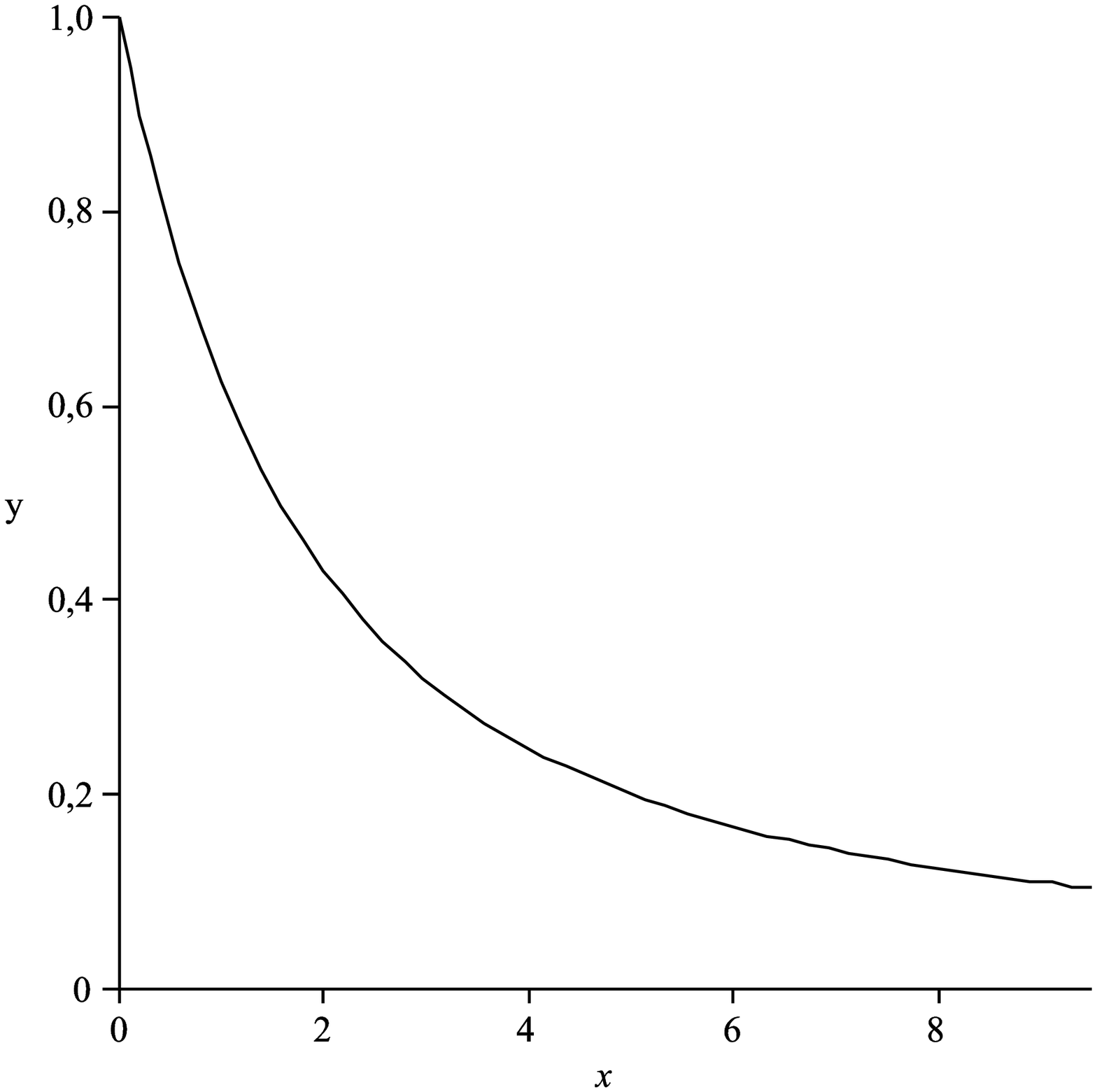}
\end{center}
 \caption{The Laplace transform of $\cA$.} 
\label{Figlap}
\end{figure}

The graph of $\psi(s)$ is shown in \refF{Figlap}.

\begin{rem} The Laplace transform of $\cA$ can also be expressed in
  terms of generalized hypergeometric functions, 
\begin{equation*}
\psi(s)= \, {}_1F_1 \Bigpar{\frac{5}{6};\frac{4}{6};\frac{s^2}{6}} -
\frac{4 s}{3 \sqrt{2 \pi}} \, {}_2F_2
\Bigpar{\frac{6}{6},\frac{8}{6};\frac{7}{6},\frac{9}{6};\frac{s^2}{6}}. 
\end{equation*}
\end{rem}

\begin{table}[ht] \newcommand{\x}{\,} %\newcommand{\x}{,}
\begin{align*}
\E\cA &= \frac{4}{3 \sqrt{2 \pi}} \x  
& 
\E\cA^2&= \frac{5}{12} \x 
& 
\E\cA^3&= \frac{64}{63 \sqrt{2 \pi}} \x  
& 
\E\cA^4&= \frac{11}{24} \x  
\end{align*}
\caption{The first four moments of $\cA$.} \label{Tabmom}
\end{table}

\section{Tail asymptotics}

Tauberian theorems
by Davies \cite{Davies} and Kasahara \cite{Kasahara} 
(see Janson \cite[Theorem 4.5]{SJ161} for a convenient version)
show that
the moment asymptotics \eqref{moms} implies the estimate
$\ln\P(\cA>x)\sim{-3x^2/2}$
for the tail of the distribution function. Thus, the following corollary is obtained.
\begin{cor}\label{tailest} $\cA$ has the tail estimate
\begin{align*}
\P(\cA>x)=\exp\bigset{-3x^2/2+o(x^2)},\qquad{x\to\infty}.
\end{align*}
\end{cor} %\vsbb
(This result can also be proved by large deviation theory; \cf{} similar results in Fill and Janson \cite{SJ197}.)

It seems difficult to obtain more precise tail
asymptotics from the moment asymptotics, but 
it is natural to make a conjecture.
\begin{conj} \label{CON1}
$\cA$ has a density function $\fa(x)$
satisfying
\begin{equation*}
\fa(x) \sim \frac{2\cdot 3^{1/6} }{\gG(2/3)}\, x^{1/3}
e^{-3x^2/2},
\qquad{x\to\infty}. 
\end{equation*}  
\end{conj}

In fact, if $\cA$ has a density with $\fa(x)\sim ax^be^{-cx^d}$ for
some constants $a,b,c,d$, then it is the only possible choice that
yields the moment asymptotics \eqref{moms}, \cf{} Janson and Louchard
\cite{SJtails}. 

Conjecture \ref{CON1} may be compared with similar results for
several Brownian areas in Janson and Louchard \cite{SJtails}, see also
Janson \cite{SJ201}. 
Note that in these result for Brownian areas, the exponent of $x$ is
always an integer (0, 1 or 2). It is therefore a small surprise that here, the exponent seems to be $1/3$, corresponding to the power $n^{1/6}$ in
\eqref{moms}.

\section{Preliminaries on point processes}\label{Sprel}

Let $\fS$ be a measurable space.
(In this paper, $\fS$ is either an interval of the real line or the
product of two such intervals.)
Although a point process $\Xi$ will be regarded as a random set 
%(or more generally, multiset) 
$\set{\xi_i}\subset\fS$, it is technically convenient to formally
define it as an integer-valued random measure
$\sum_i\delta_{\xi_i}$.
Hence, $\Xi(A)$ denotes the number of points $\xi_i$ that belong to a (measurable)
subset $A\subseteq\fS$.
Also,
$x\in\Xi$ is equivalent to $\Xi(\set{x})>0$.
See further \eg{} Kallenberg \cite{Kallenberg}.

A \emph{Poisson process with intensity $\dd\mu$}, where $\dd\mu$ is a
measure on $\fS$, is a point process $\Xi$ such that 
$\Xi(A)$ has a Poisson distribution with mean $\mu(A)$
for every measurable
$A\subseteq\fS$, 
and $\Xi(A_1),\dots,\Xi(A_k)$ are independent for every family
$A_1,\dots,A_k$ of disjoint measurable sets.
Lemma \ref{Lpo} is a standard formula for Laplace functionals, see for instance \cite[Lemma 12.2(i)]{Kallenberg}. 

\begin{lem} \label{Lpo}
If\/ $\Xi$ is a Poisson process with intensity $\dd\mu$ on a set $\fS$,
and $f:\fS\to[0,\infty)$ is a measurable function, then
\begin{equation*}
\E \exp\Bigset{-\sum_{\xi\in\Xi} f(\xi)} =
\exp\Bigset{-\int_{\fS}\bigpar{1-e^{-f(x)}}\dd\mu(x)}.
\end{equation*}
\end{lem} %\vsbb
Lemma \ref{Lgamma}, on the other hand, is more of a digression. The
result follows 
from a standard Gamma integral by integration by parts.
(The result can also be written as $2\Gamma\bigpar{1/2}\gl\qq$.)

\begin{lem}
  \label{Lgamma}
If $\gl>0$, then
\begin{equation*}
  \intoo\bigpar{1-\expe{-\gl x}} x\qqcw\dd x
=2\sqrt{\pi\gl}.
\end{equation*}
\end{lem} %\vsbb

\section{Proof of \refT{T1}}\label{Spf}

The set \set{t:B(t)=0} is \as{} closed and unbounded, so its
complement \set{t:B(t)\neq0} is an infinite union of finite open
intervals, denoted by $\inu=(g\nnu,d\nnu)$, $\nu=1,2,\dots$,
in some order. (The intervals cannot be ordered by appearance, 
since there is \as{} an infinite number of them in, say,
\oi. Fortunately, the order does not matter.)
The restrictions of $B(t)$ to these intervals are called the
\emph{excursions} of $B(t)$. Let $\enu$ be the excursion during $\inu$.

The local time $L(t)$ is constant during each excursion.
Let $\tau\nnu$ be the local time during $\enu$
and let $\ell\nnu\=d\nnu-g\nnu$ be the length of $\enu$.
It is well-known, see Revuz and Yor \cite[Chapter XII]{RY}, that the collection of
pairs $\set{(\tau\nnu,\ell\nnu)}\nuioo$
forms a Poisson process in $[0,\infty)\times(0,\infty)$ with intensity
  \begin{equation*}
\dd\gL=(2\pi\ell^3)\qqw\dd \tau\dd \ell.	
  \end{equation*}
Note also that, \as, if the excursion $\exc_{\nu_1}$ comes before 
$\exc_{\nu_2}$, then 
$\tau_{\nu_1}<\tau_{\nu_2}$.

Next, consider a Poisson process $\set{T_i}\iioo$ on $\ooo$ with
intensity $\gl\dd t$, independent of $\set{B(t)}$. Assume that the
points are ordered with $0<T_1<T_2<\dotsm$. Then $T_1$,
$T_2-T_1$, \dots\ are \iid{} 
$\Exp(\gl)$ random variables
with density function $\gl\expe{-\gl t}$.
Furthermore, $T_1$ is independent of $\set{B(t)}$ and thus of
$\set{\cA(T)}$. It follows from \eqref{a2} that
$  \cA(T_1)\eqd T_1\qqc\cA$ and consequently
\begin{equation}\label{ea1}
\E \expe{-\ga\cA(T_1)}
=
\E \expe{-\ga T_1\qqc\cA}
=
\E \psi\bigpar{\ga T_1\qqc}
=
\gl\intoo e^{-\gl s} \psi\bigpar{\ga s\qqc}\dd s.
\end{equation}

The times $T_i$ are called \emph{marks}, and
an excursion is called \emph{marked} if it contains at least one
of the marks $T_i$.
The marks \set{T_i} are placed by first constructing \set{B(t)}
and then adding marks according to independent Poisson processes with
intensities $\gl\dd t$ in each excursion.
Thus, given the excursions \set{\enu}, each excursion  $\enu$
is marked with probability $1-\expe{-\gl\ell\nnu}$, independently of
the other excursions. The Poisson process $\Xi\=\set{(\tau\nnu,\ell\nnu)}$ 
defined by the excursions can be written as the union $\Xi'\cup\Xi''$, where
\begin{align*}
  \Xi'&\=\bigset{(\tau\nnu,\ell\nnu):\enu\text{ is unmarked}},
\\
  \Xi''&\=\bigset{(\tau\nnu,\ell\nnu):\enu\text{ is marked}}.
\end{align*}
By the general independence properties of Poisson processes, 
$\Xi'$ and $\Xi''$ are \emph{independent} Poisson processes
with intensities 
\begin{align}\label{gL'}
\dd\gL'&\=
\expe{-\gl\ell}\dd\gL
=(2\pi)\qqw\ell\qqcw\expe{-\gl\ell}\dd \tau\dd \ell  \\
\intertext{and}\label{gL''}
\dd\gL''&\=
\bigpar{1-\expe{-\gl\ell}}\dd\gL
=(2\pi)\qqw\ell\qqcw\bigpar{1-\expe{-\gl\ell}}\dd \tau\dd \ell, 
\end{align}
respectively.
In particular, if the lengths are ignored, the local times of the marked
excursions form a Poisson process $\tXi$ on $(0,\infty)$ with
intensity
\begin{equation*}
\int_{\ell=0}^\infty \bigpar{1-\expe{-\gl\ell}}\dd\gL
=\tgl\dd \tau,
\end{equation*}
where, using Lemma \ref{Lgamma},
\begin{equation}\label{tgl}
\tgl=\intoo(2\pi)\qqw\ell\qqcw\bigpar{1-\expe{-\gl\ell}}\dd \ell
=\sqrt{2\gl}.
\end{equation}

Due to the fact that $B(T_1)\neq0$ \as, there exists a unique excursion $\exc\nnux$
that contains the first mark $T_1$, \ie, $T_1\in I\nnux$.
Let $\tx\=L(T_1)=\tau\nnux$ be the local time at $T_1$ (and thus during
$\enux$).
Since $\enux$ is the first marked excursion, its local time $\tx$ is
the first of the points in the Poisson process $\tXi$ and hence
\begin{equation}\label{tx}
  \tx\sim\Exp\bigpar{\sqrt{2\gl}}.
\end{equation}

The restriction of $B(t)$ to the interval $[0,T_1]$  consists of
all excursions $\enu$ with local time $\tau\nnu<\tau\nnux=\tx$ and the part of
$\enux$ on $(g\nnux,T_1)$, plus the set
\begin{equation*}
[0,T_1]\setminus\bigcup_\nu \inu = \set{t\le T_1: B(t)=0}
\end{equation*}
which \as{} has measure 0 and thus may be ignored.
Consequently, since $L(t)=\tau\nnu$ on $\inu$,
\begin{equation*}
  \begin{split}
\cA(T_1)
&\=\int_0^{T_1} L(t)\dd t
=\sum_{\nu:\tau\nnu<\tau\nnux} \int_{\inu}L(t)\dd t
+ \int_{g\nnux}^{T_1}L(t)\dd t	
\\&\phantom{:}
=\sum_{\nu:\tau\nnu<\tx} \tau\nnu\ell\nnu + \tx(T_1-g\nnux) :=\cA'+\cA''.
  \end{split}
\end{equation*}

The sum defined as $\cA'=\sum_{\nu:\tau\nnu<\tx} \tau\nnu\ell\nnu$
only contains terms for unmarked excursions $\enu$. Thus
\begin{equation*}
\cA'
  =\sum_{(\tau\nnu,\ell\nnu)\in\Xi':\;\tau\nnu<\tx} \tau\nnu\ell\nnu.
\end{equation*}
Recall that $\zeta$ is determined by $\Xi''$ (as the smallest $\tau$
with $(\tau,\ell)\in\Xi''$ for some $\ell$)
and that $\Xi'$ and $\Xi''$ are independent. Hence, $\Xi'$ and $\tx$
are independent.
It follows from \refL{Lpo}, with 
$\fS=(0,\tx)\times(0,\infty)$ and 
$f((\tau,\ell))=\ga\tau\ell$,
that
\begin{equation*}
\E\bigpar{e^{-\ga\cA'}\bigm|\tx} = \exp\Bigset{-\int_{\tau=0}^\tx
  \int_{\ell=0}^\infty \bigpar{1-e^{-\ga \tau\ell}}\dd\gL'(\tau,\ell)}. 
\end{equation*}
By \eqref{gL'} and Lemma \ref{Lgamma},
\begin{equation*}
\begin{split}
\int_{\tau=0}^\tx \int_{\ell=0}^\infty \bigpar{1 & -e^{-\ga  \tau\ell}} 
\dd\gL'(\tau,\ell) \\ 
& = \int_{\tau=0}^\tx \int_{\ell=0}^\infty \bigpar{1-e^{-\ga \tau\ell}}
(2\pi)\qqw\ell\qqcw\expe{-\gl\ell}\dd \ell\dd \tau  \\
& = (2\pi)\qqw\int_{\tau=0}^\tx \int_{\ell=0}^\infty
\bigpar{\expe{-\gl\ell}-e^{-(\gl+\ga \tau)\ell}} \ell\qqcw\dd \ell\dd
\tau \\ 
& = \int_{\tau=0}^\tx \sqrt2\lrpar{\sqrt{\gl+\ga\tau}-\sqrt{\gl}}\dd\tau
\\
& = \frac{2\sqrt2}{3\ga} \lrpar{(\gl+\ga\tx)\qqc-\gl\qqc} -\sqrt{2\gl}\,\tx,
\end{split}
\end{equation*}
and it follows that
\begin{equation}\label{ea'}
\E\bigpar{e^{-\ga\cA'}\bigm|\tx} = \exp\Bigset{\sqrt{2\gl}\,\tx
  -\frac{2\sqrt2}{3\ga} \lrpar{(\gl+\ga\tx)\qqc-\gl\qqc}}. 
\end{equation}

Now consider $\cA''=\tx(T_1-g\nnux)$. Note that $T_1-g\nnux$ is the
location (relative to the left endpoint of the excursion) of the first
mark in the first marked excursion. Since $\Xi$ is a Poisson
process with intensity independent of $\tau$, the location $T_1-g\nnux$ is
independent of the local time $\tx$ of the first marked excursion.
Further, the joint distribution of $(\ell\nnux,T_1-g\nnux)$ has density
\begin{equation*}
(\tgl)\qw \gl \expe{-\gl y}(2\pi)\qqw\ell\qqcw\dd \ell \dd y ,
\qquad 0<y<\ell<\infty,
\end{equation*}
where the normalization constant $\tgl$ is given by \eqref{tgl}.
Consequently,
\begin{equation}\label{ea''}
  \begin{split}
\E\bigpar{e^{-\ga\cA''}\bigm|\tx}
&=
\E\bigpar{e^{-\ga\tx(T_1-g\nnux)}\bigm|\tx}
\\&
=
\int_{y=0}^\infty \int_{\ell=y}^\infty 
e^{-\ga \tx y}
(\tgl)\qw \gl \expe{-\gl y}(2\pi)\qqw\ell\qqcw\dd \ell \dd y 
\\&
=\pi\qqw\gl\qq	\int_{y=0}^\infty 
 \expe{-(\gl+\ga\tx) y} y\qqw \dd y 
\\&
=\gl\qq	(\gl+\ga\tx)\qqw.
  \end{split}
\end{equation}

Again, since $\Xi'$ and $\Xi''$ are independent, $\cA'$ and
$\cA''$ are conditionally independent given $\tx$. Thus, equation \eqref{ea'} and \eqref{ea''} yield
\begin{equation*}
  \begin{split}
\E\bigpar{e^{-\ga\cA(T_1)}\bigm|\tx} &=
\E\bigpar{e^{-\ga\cA'}\bigm|\tx} \E\bigpar{e^{-\ga\cA''}\bigm|\tx} \\ 
& = \parfrac{\gl}{\gl+\ga\tx}\qq \exp \Bigset{\sqrt{2\gl}\,\tx
-\frac{2\sqrt2}{3\ga}
\lrpar{(\gl+\ga\tx)\qqc-\gl\qqc}}.
  \end{split}
\end{equation*}
By \eqref{tx}, $\tx$ has the density $\sqrt{2\gl} e^{-\sqrt{2\gl}x}$,
$x>0$, and it follows that 
\begin{equation*}
\E e^{-\ga\cA(T_1)} 
= \gl\sqrt2\intoo(\gl+\ga x)\qqw \exp\Bigset{-\frac{2\sqrt2}{3\ga}
  \lrpar{(\gl+\ga x)\qqc-\gl\qqc}} \dd x. 
\end{equation*}
Finally, the substitution 
\begin{equation*}
\frac{2\sqrt2}{3\ga \gl} \lrpar{(\gl+\ga x)\qqc-\gl\qqc} \mapsto s
\end{equation*}
provides the slightly simpler formula
\begin{equation*}
\E e^{-\ga\cA(T_1)} = \gl \intoo \Bigpar{1+\frac{3 \ga s}{2 \sqrt{2
	  \gl}}}^{-2/3}  e^{-\gl s} \dd s. 
\end{equation*} 
The result now follows by a comparison with \eqref{ea1}.
\qed

\newcommand\webcite[1]{\hfil
   \penalty0\texttt{\def~{{\tiny$\sim$}}#1}\hfill\hfill}
\newcommand\arxiv[1]{\webcite{arXiv:#1.}}

\end{document}